\documentclass[10pt,oneside,a4paper,reqno]{amsart}  

\usepackage{amssymb,graphicx,xcolor}
\usepackage[colorlinks=true, linkcolor=magenta, citecolor=cyan, urlcolor=blue]{hyperref}
\usepackage{bbm}

\usepackage{graphicx}
\usepackage{mathrsfs}

\numberwithin{equation}{section}\newtheorem{theorem}{Theorem}[section]
\newtheorem{lemma}[theorem]{Lemma}

\theoremstyle{remark}
\newtheorem{remark}{Remark}[section]
\theoremstyle{definition}

\newcommand{\p}{\widetilde{p}}

\newcommand{\Rn}{\mathbb{R}^{n}}
\newcommand{\Rpiu}{\mathbb{R}^{+}}

\title{Singular integrals with angular integrability}
\date{}    

\author{Federico Cacciafesta}
\address{Federico Cacciafesta: 
SAPIENZA --- Universit\`a di Roma,
Dipartimento di Matematica, 
Piazzale A.~Moro 2, I-00185 Roma, Italy}
\email{cacciafe@mat.uniroma1.it}

\author{Renato Luc\`a}
\address{Renato Luc\`a: Instituto de Ciencias Matem\'aticas CSIC-UAM-UC3M-UCM, Madrid, 28049, Spain.}
\email{renato.luca@icmat.es}

\subjclass[2010]{42B37, 42B20}
\keywords{Singular integrals. Angular integrability}

\begin{document}

 \begin{abstract}
In this note we prove a class of sharp inequalities for singular integral operators in weighted Lebesgue spaces with angular integrability.

\end{abstract}

\maketitle

\section{Introduction}

We consider singular integral operators
\begin{equation}\label{singint}
T f (x):=
\mathrm{P.V.}
\int_{\mathbb{R}^{n}} f(x-y) K(y) \ dy,
\end{equation}
where the kernel $K$ satisfies the following conditions
$$
|y|^{n} |K| \leq C,
\quad
|y|^{n+1} | \nabla K | \leq C,
\quad
| \widehat{K} | \leq C.
$$
Here $C >0 $ is a constant and $\, \widehat{} \,$ denotes the Fourier transform.
The main example we have in mind is the directional Riesz transform, which corresponds to the choice 
$K(y) := |y|^{-(n+1)} y \cdot \theta$, $\theta \in \mathbb{S}^{n-1}$.

The study of the boundedness of these operators in weighted Lebesgue spaces $L^p(w(x)dx)$,  for $1 < p <\infty$ and $0 < w \in L^{1}_{loc}(\Rn)$, is a classical problem in harmonic analysis: in particular, Stein \cite{Stein} proved it for the (sharp range of) homogeneous weights $w(x)=|x|^{\alpha p}$, $-n/p < \alpha < n- n/p$. The result was later extended by Coifman and Fefferman \cite{coiffeff} to any $A_p$ weight.
%

While the weighted $L^{p}$-theory has been extensively studied, less is known in the case of Lebesgue norms with different integrability in the radial and angular directions, namely

\begin{equation}\label{MixedNorms}
\begin{array}{lcl}
\displaystyle
  \|f\|_{L^{p}_{|x|}L^{\p}_{\theta}}& := &
\displaystyle
 \left(
    \int_{0}^{+\infty}
    \|f(\rho \, \cdot \, )\|^{p}_{L^{\p}(\mathbb{S}^{n-1})}
    \rho^{n-1}d \rho
 \right)^{\frac1p}. 
\end{array}
\end{equation}

These mixed radial-angular spaces have been successfully used in recent years to improve several results in the framework of 
partial differential equations; see e.g. \cite{cacdan, FangWang08-a, Ren, MachiharaNakamuraNakanishi05-a, Rogers, Sterbenz05-a, TaoSphere}. 
Notice that when $p=\p$ the norms reduce to the usual $L^{p}$ norms. 
Notice also that, neglecting the constants, they are increasing in $\p$, and that they behave as the $L^{p}$ norms under homogeneous rescaling, namely
$f(\cdot) \to f(\lambda \, \cdot \, )$, $\lambda > 0$.

In a recent paper A. C\'{o}rdoba \cite{Cordoba} proved, among the other things, 
the 
$L^{p}_{|x|}L^{2}_{\theta}$ boundedness for operators of the form (\ref{singint}).
Here we give an extension of this result to the weighted setting.

\begin{theorem}\label{OurSteinWeiss}
Let $n \geq 2$, $1 < p < \infty$, $1 < \p < \infty$ and $-n/p < \alpha < n- n/p$. 
Then
\begin{equation}\label{eq:OurSteinWeiss}
\| |x|^{\alpha} T \phi \|_{L^{p}_{|x|}L^{\p}_{\theta}} 
\leq C \| |x|^{\alpha} \phi \|_{L^{p}_{|x|}L^{\p}_{\theta}},
\end{equation}  
where $C$ is a constant depending only on $\alpha, p, \p, n$.
\end{theorem}

Let us point out that the case $\alpha =0$, $1 < \p < \infty$ in inequality (\ref{eq:OurSteinWeiss}) may be deduced by the application of C\'{o}rdoba's argument; see \cite[Theorem 2.6]{DLNS}. 
Therefore the novelty of Theorem \ref{OurSteinWeiss} is in that it covers all the possible homogeneous weights of the kind $|x|^{\alpha p}$.

\begin{remark}
Condition $-n/p < \alpha$ turns to be necessary by testing the inequality on functions $\phi$ such that $\phi=0$ and $T\phi>0$ in a neighborhood of the origin. On the other hand, condition $\alpha < n - n/p$ turns to be necessary for the same reason by considering the dual inequality.
\end{remark}

\begin{remark}
One of the estimates (\ref{eq:OurSteinWeiss}) has been used in \cite[Theorem 1.5]{DLNS} to deduce informations
about the regularity of weak solutions of the $3d$ Navier--Stokes problem with initial velocities
satisfying good angular integrability properties.
\end{remark}

We write $A \lesssim B$ if $A \leq C B$ with a constant $C$ depending only on $\alpha, p, \p, n$. We write $A \simeq B$ if both $A \lesssim B$ and $B \lesssim A$.

\section{Proof}
\setcounter{equation}{0}

We know by \cite[Theorem 2.6]{DLNS} that inequality (\ref{eq:OurSteinWeiss}) is true in the case $\alpha =0$, that is
\begin{equation}\label{Eq:NONweightedSW}
\| T g \|_{L^{p}_{|x|}L^{\p}_{\theta}} \lesssim \| g \|_{L^{p}_{|x|}L^{\p}_{\theta}}.
\end{equation}

Following Stein \cite{Stein}, we now show that the weighted case (\ref{eq:OurSteinWeiss}) can be then deduced by the unweighted one. The next lemma represents the core of the proof.

\begin{lemma}\label{SteinMainLemma}
Let $n \geq 2$, $1 < p < \infty$, $1\leq \p \leq \infty$, $-n/p < \alpha < n - n/p$ and
\begin{equation}\label{eq:SteinWeissKernel}
F(x,y) := \frac{|1- (|x|/|y|)^{\alpha}|}{|x-y|^{n}},
\end{equation}
then
\begin{equation}\label{Eq:SWInequalityGeneralized}
\bigg\| \int_{\Rn} F(x,y) \phi(y) dy  \bigg\|_{L^{p}_{|x|}L^{\p}_{\theta}}
\leq C 
\| \phi \|_{L^{p}_{|x|}L^{\p}_{\theta}}.
\end{equation}
\end{lemma}

Assume indeed this has been proved and
first apply inequality \eqref{Eq:NONweightedSW} with the choice $g:= |\cdot|^{\alpha} f$ to have
\begin{equation}\label{SWQuasiI}
\| T (|x|^{\alpha}f) \|_{L^{p}_{|x|}L^{\p}_{\theta}} \lesssim \| |x|^{\alpha} f \|_{L^{p}_{|x|}L^{\p}_{\theta}}.
\end{equation}
Then notice that  
\begin{equation}\label{SWQuasiII}
\begin{split}
& \big| T(|x|^{\alpha} f)  -   |x|^{\alpha} T f \big| 
 \leq
 \int_{\Rn}  \big| K(x-y) ( |y|^{\alpha} - |x|^{\alpha}) f(y) \big|  dy  
\\ 
& \lesssim
 \int_{\Rn}    \frac{| |y|^{\alpha} - |x|^{\alpha} |}{|x-y|^{n}}  |f(y)|  dy  
=
 \int_{\Rn} \frac{|1 - (|x|/|y|)^{\alpha}|}{|x-y|^{n}} |y|^{\alpha}  |f(y)|  dy  ,
\end{split}
\end{equation}
so that by using Lemma \ref{SteinMainLemma} with $\phi = |\cdot| ^{\alpha} f$ we obtain
\begin{equation}
\|  T(|x|^{\alpha} f)  -   |x|^{\alpha} T f  \|_{L^{p}_{|x|}L^{\p}_{\theta}}
\lesssim \| |x|^{\alpha} f\|_{L^{p}_{|x|}L^{\p}_{\theta}} .
\end{equation}
Then, the desired estimate   
\eqref{eq:OurSteinWeiss} follows by (\ref{SWQuasiI}, \ref{SWQuasiII}) and triangle inequality.
Thus it only remains to prove Lemma \ref{SteinMainLemma}.

The idea is to use a change of variables which resembles the standard polar coordinates. 
In this variant the integration over the sphere is replaced by integration over the special orthogonal group $SO(n)$ and the
radial integration is replaced by integration over the multiplicative group of the positive real numbers. 
This method works efficiently 
when homogeneous power weights are involved; see e.g. \cite{DL, DenapoliDrelichmanDuran11-a}.

\subsubsection*{Proof of Lemma \ref{SteinMainLemma}}

By using the isomorphism 
\begin{equation*}
  \mathbb{S}^{n-1}\simeq SO(n)/SO(n-1)
\end{equation*}
we can rewrite integrals on $\mathbb{S}^{n-1}$ as follows
\begin{equation*}
  \int_{\mathbb{S}^{n-1}}g(y)dS(y)
  \simeq
  \int_{SO(n)}g(Ae)dA,\qquad n\ge2,
\end{equation*}
where $dA$ is the left Haar measure on
$SO(n)$, and $e\in\mathbb{S}^{n-1}$ is a fixed 
unit vector. Thus, via polar coordinates, 
a generic integral can be rewritten as
\begin{equation*}
\begin{split}
  \int_{\Rn} F(x,y) \phi(x) dy
   &=
  \int_{0}^{\infty}
  \int_{\mathbb{S}^{n-1}}
     F(x, \rho \omega)\phi(\rho \omega)dS_{\omega}\rho^{n-1}d \rho
    \\
  &\simeq
  \int_{0}^{\infty}
  \int_{SO(n)}F(x, \rho Be)\phi(\rho Be)dB\rho^{n-1}d \rho.
\end{split}
\end{equation*}
Hence the $L^{\p}_{\theta}$ norm can be written as
\begin{equation*}
\begin{split}
  & \bigg\| \int_{\Rn} F(|x|\theta,y) \phi(y) dy\bigg\|_{L^{\p}_{\theta}(\mathbb{S}^{n-1})}
  \simeq
  \bigg\| \int_{\Rn} F(|x| Ae,y) \phi(y) dy \bigg\|_{L^{\p}_{A}(SO(n))}
    \\
  &\le 
  \int_{0}^{\infty}
  \bigg\|
    \int_{SO(n)}F(|x|Ae, \rho Be)\phi(\rho Be)dB
  \bigg\|_{L^{\p}_{A}(SO(n))}
  \rho^{n-1}d \rho
\end{split} 
\end{equation*}
where $e$ is any fixed unit vector.
We choose $F$ as in \eqref{eq:SteinWeissKernel}
and we change variables
$B\to AB^{-1}$ in the inner integral.
By the invariance of the measure 
this is equivalent to
\begin{equation*}
\begin{split}
  & =\int_{0}^{\infty}
  \bigg\|
    \int_{SO(n)}\frac{|1 - (|x|/\rho)^{\beta}|}{|AB^{-1}(|x|Be-\rho e)|^n}\phi(\rho AB^{-1}e)dB
  \bigg\|_{L^{\p}_{A}(SO(n))}
  \rho^{n-1}d \rho
\\ 
&
=\int_{0}^{\infty}
  \bigg\|
    \int_{SO(n)}\frac{|1 - (|x|/\rho)^{\beta}|}{||x|Be-\rho e|^n}\phi(\rho AB^{-1}e)dB
  \bigg\|_{L^{\p}_{A}(SO(n))}
  \rho^{n-1}d \rho.
\end{split}
\end{equation*}
Notice that the integral
\begin{equation*}
    \int_{SO(n)}\frac{|1 - (|x|/\rho)^{\beta}|}{||x|Be-\rho e|^n}  \  |\phi(\rho AB^{-1}e)|dB
  =
  G*\phi(A)
\end{equation*}
is a convolution on $SO(n)$ of the functions
\begin{equation*}
  G(A)=\frac{|1 - (|x|/\rho)^{\beta}|}{||x|Ae-\rho e|^n},
  \qquad
  H(A)=|\phi(\rho Ae)|.
\end{equation*}
We can thus apply Young's inequality on $SO(n)$
(see for instance  
\cite[Theorem 1.2.12]{Grafakos08-a}) to obtain, for any
$1 \leq \p \leq \infty$,
the estimate
\begin{equation}\label{eq:firstest}
\begin{split}
 & \left\|\int_{\Rn}F(|x|\theta,y)\phi(y) dy\right\|_{L^{\p}_{\theta}(\mathbb{S}^{n-1})}
 \!\!\!
 \lesssim
  \int_{0}^{\infty}
  \left\| \frac{|1 - (|x|/\rho)^{\beta}|}{|(|x|e-\rho \theta)|^n} \right\|_{L^{1}_{\theta}(\mathbb{S}^{n-1})}
  \!\!\!\!\!\!\!\!\!\!\!\!\!
  \|\phi(\rho\theta)\|_{L^{\p}_{\theta}(\mathbb{S}^{n-1})}
  \rho^{n-1}d \rho
  \\ 
 & = 
  \int_{0}^{\infty}
  \bigg\| \frac{|1 - (|x|/\rho)^{\beta}|}{|(\frac{|x|}{\rho} -  \theta)|^n} \bigg\|_{L^{1}_{\theta}(\mathbb{S}^{n-1})}
  \!\!\!\!\!\!\!\!\!\!\!\!    \|\phi(\rho\theta)\|_{L^{\p}_{\theta}(\mathbb{S}^{n-1})}
 \frac{d \rho}{\rho}
  \end{split}
\end{equation}
where we switched back to the coordinates of $\mathbb{S}^{n-1}$.
Then we notice
\begin{equation}\label{est}
\begin{split}
\left \| \int_{\Rn} F(x,y) \phi(y)  dy  \right\|_{L^{p}_{|x|}L^{\p}_{\theta}} 
& 
\lesssim
\\ 
 &
\!\!\!\!\!\!\!\!\!\!\!\!\!\!\!\!\!\!\!\!\!\!\!\!\!\!\!\!\!\!\!\!\!\!\!\!\!\!\!\!\!
 \lesssim 
 \Bigg\| \bigg\| \int_{\Rn} F(|x| \theta,y) \phi(y) dy  \bigg\|_{L^{\p}_{\theta}(\mathbb{S}^{n-1})} \!\!\!\!\!\!\!\!\!\!    dS_{\theta} |x|^{\frac np}  \Bigg\|_{L^{p}\left(\Rpiu(\cdot), d|x|/|x|\right)}
\end{split}
\end{equation}
where $\Rpiu(\cdot)$ is the moltiplicative group of positive real numbers equipped with its Haar measure $d\rho/\rho $.
Using \eqref{eq:firstest} allows to estimate \eqref{est} with
\begin{equation}\nonumber
\lesssim 
  \Bigg\|  
  \int_{0}^{\infty}
\left( \frac{|x|}{\rho} \right)^{\frac np}  
\bigg\|
\frac{|1 - (|x|/\rho)^{\beta}|}{\big| \big(\frac{|x|}{\rho} e - \theta \big) \big|^n} \bigg\|_{L^{1}_{\theta}(\mathbb{S}^{n-1})}
\!\!\!\!\!\!\!\!\!\!\!
\rho^{\frac np}
\|\phi(\rho\theta)\|_{L^{\p}_{\theta}(\mathbb{S}^{n-1})}
 \frac{d \rho}{\rho}
\Bigg\|_{L^{p}\left(\Rpiu(\cdot), d|x|/|x|\right)}.
\end{equation}
Notice that the inner term is a convolution on 
$\Rpiu (\cdot)$ of the functions
\begin{equation}\nonumber
g(\rho) = \rho^{\frac np}  
\left\| \frac{|1 - \rho^{\beta}|}{| (\rho e - \theta )|^n} \right\|_{L^{1}_{\theta}(\mathbb{S}^{n-1})},
\qquad
h(\rho) =\rho^{\frac np}
\|\phi(\rho\theta)\|_{L^{\p}_{\theta}(\mathbb{S}^{n-1})}.
\end{equation}
Thus we can apply again Young's inequality to estimate further (\ref{est}) with
\begin{equation}\label{eq:SecondYoung}
 \lesssim
\big\|  \|  g(\rho) \|_{L^{1}_{\theta}(\mathbb{S}^{n-1})} \big\|_{L^{1}(\Rpiu (\cdot), d \rho / \rho)} 
\big\|  \|  h(\rho) \|_{L^{\p}_{\theta}(\mathbb{S}^{n-1})} \big\|_{L^{p}(\Rpiu (\cdot), d \rho / \rho)},
\end{equation}
for all $1 \leq p \leq \infty$.
Once we have noticed that
\begin{equation}\nonumber
\big\|  \|  h(\rho) \|_{L^{\p}_{\theta}(\mathbb{S}^{n-1})} \big\|_{L^{p}(\Rpiu (\cdot), d \rho / \rho)}
= \| \phi \|_{L^{p}_{|x|}L^{\p}_{\theta}},
\end{equation}
the concluding step of the proof of the Lemma is represented by showing that the first term of \eqref{eq:SecondYoung} is bounded.
We split the integral 
\begin{equation}\nonumber
\int_{0}^{+\infty} \rho^{\frac np}  \int_{\mathbb{S}^{n-1}} \frac{|1 - \rho^{\beta}|}{|\rho e -\theta|^{n}} \  dS_{\theta}  \  \frac{d \rho}{\rho}
= \int_{0}^{\frac 12} (\cdot)
+
\int_{\frac 12}^{2} (\cdot)
+
\int_{2}^{+\infty} (\cdot)
=: I + II + III
\end{equation} 
and we bound separately the three terms.

If $0< \rho < 1/2$, then $| \rho e - \theta | \geq 1/2 $. Thus, since $|1 - \rho^{\beta}| < 1 + \rho^{\beta}$ 
\begin{equation}\nonumber
I \lesssim \int_{0}^{\frac 12} ( \rho^{\frac np -1} + \rho^{\frac np - 1 + \beta} ) d \rho < \infty
\end{equation}
provided that $p < \infty$ and $\beta > - n/p$.

If $1/2 \leq \rho \leq 2$, we notice that $|1 - \rho^{\beta}| \lesssim | 1-\rho |$ and we use that
(see for instance Lemma 2.1 in \cite{DL})
\begin{equation}\nonumber
\int_{\mathbb{S}^{n-1}} |\rho e - \theta|^{-n} dS_{\theta}
 \simeq \frac{1}{|1- \rho |}
\end{equation} 
to bound
\begin{equation}
II \simeq
\int_{\frac 12}^{2} \rho^{\frac np -1} \frac{|1 - \rho^{\beta}|}{|1- \rho|} d \rho 
\lesssim
\int_{\frac 12}^{2} \rho^{\frac np -1}  d \rho
<
\infty.
\end{equation}

Finally, if $2 < \rho < + \infty$, then $| \rho e - \theta | \geq | \rho | - | \theta | \geq | \rho | /2 $. Thus, since $|1 - \rho^{\beta}| < 1 + \rho^{\beta}$ 
\begin{equation}\nonumber
III \lesssim \int_{2}^{+ \infty} ( \rho^{\frac np -1 -n} + \rho^{\frac np - 1 + \beta - n} ) d \rho < \infty
\end{equation}
provided that $p > 1$ and $\beta < n - n/p$,
that concludes the proof. 

\hfill $\Box$

{\bf Acknowledgements.} The first author is supported by the FIRB 2012 \lq Dispersive dynamics, Fourier analysis and variational methods'. Part of the work was done during the first author's visit at MSRI, Berkeley (CA), within the 
program \lq New Challenges in PDE: Deterministic Dynamics and Randomness in High and Infinite Dimensional Systems', which he acknowledges for the wonderful working conditions. 

The second author is supported by the
ERC grant 277778 and MINECO grant SEV-2011-0087 (Spain).

\end{document}